\documentclass[10pt,twoside]{article}
\usepackage{graphicx}
\usepackage{amsmath,amssymb}
\usepackage{Latex-document}

\title{\bf Rhythms of the Nervous System: \vskip -2mm
 Mathematical Themes and Variations\vskip 6mm}

\author{Nancy Kopell\thanks{Boston University, Department of Mathematics and Center for BioDynamics,
Boston, MA 02215, USA. E-mail: nk@bu.edu}\vspace*{-0.5cm}}
\date{\vspace{-8mm}}

\begin{document}
\maketitle

\thispagestyle{first} \setcounter{page}{805}

\begin{abstract}\vskip 3mm
The nervous system displays a variety of rhythms in both waking and sleep.
These rhythms have been closely associated with different behavioral and
cognitive states, but it is still unknown how the nervous system makes use
of these rhythms to perform functionally important tasks.  To address
those questions, it is first useful to understood in a mechanistic way the
origin of the rhythms, their interactions,  the signals which create the
transitions among rhythms, and the ways in which  rhythms filter the
signals to a network of neurons.

This talk discusses how dynamical systems have been used to investigate
the origin, properties and interactions of rhythms in the nervous system.
It focuses on how the underlying physiology of the cells and synapses of
the networks shape the dynamics of the network in different contexts,
allowing the variety of dynamical behaviors to be displayed by the same
network. The work is presented using a series of related case studies on
different rhythms.  These case studies are chosen to highlight
mathematical issues, and suggest further mathematical work to be done.
The topics include: different roles of excitation and inhibition in
creating synchronous assemblies of cells, different kinds of building
blocks for neural oscillations, and transitions among rhythms.  The
mathematical issues include reduction of large networks to low dimensional
maps, role of noise, global bifurcations, use of probabilistic
formulations.

\vskip 4.5mm

\end{abstract}

\vskip 12mm

\section{Introduction} \setzero \vskip-5mm \hspace{5mm }

The nervous system creates many different rhythms, each associated with a
range of behaviors and cognitive states. The rhythms were first discovered
from scalp recordings of humans, and the names by which they are known still
come mainly from the electroencephalograph (EEG) literature, which pays
attention to the frequency and behavioral context of those rhythms, but not
to their mechanistic origins. The rhythmic patterns include the alpha (9-11
Hz), beta (12-30 Hz), gamma (30-80 Hz), theta (4-8 Hz), delta (2-4) Hz and
slow wave (.5-2 Hz) rhythms. The boundaries of these ranges are rough. More
will be said about the circumstances in which some of these rhythms are
displayed.

It is now possible to get far more information about the mechanisms behind
the dynamics of the nervous system from other techniques, including
electrophysiology. The revolutions in experimental techniques, data
acquisition and analysis, and fast computation have opened up a broad and
deep avenue for mathematical analysis. The general question addressed by
those interested in rhythms is: how does the brain make use of these rhythms
in sensory processing, sensory-motor coordination and cognition? The
mathematical strategy, to be discussed below, is to investigate the
``dynamical structure'' of the different rhythms to get clues to function.
Most of this talk is about dynamical structure, and the mathematical issues
surrounding its investigation. I'll return at the end to issues of function.

\section{Neuromath}\setzero \vskip-5mm \hspace{5mm }

The mathematical framework for the study of brain dynamics are the
Hodgkin-Huxley (HH) equations. These are partial differential equations
describing the propagation of action potentials in neurons (cells of the
nervous system). The equations, which play the same role in neural dynamics
that Navier-Stokes does in fluid dynamics, are an elaborate analogy to a
distributed electrical circuit. The central equation,

\[
Cv^{\prime }=\sum I_{ion}+\triangledown ^{2}v+\sum I_{syn}+I_{ext}
\]
describes conservation of current across a piece of a cell membrane; $v$ is
the cross-membrane voltage and the left hand side is the capacitive current.
The first sum on the right hand side represents the intrinsically generated
ionic currents across the membrane. The term $\triangledown ^{2}v$
represents the spatial diffusion, and $I_{ext}$ the current fed into the
cell. $\sum I_{syn}$ represents the currents introduced by coupling from
other cells. Thus, these equations can also be used to model networks of
interacting neurons, the focus of this talk.

Each of the intrinsic currents $I_{ion}$ is described by Ohm's law: $I_{ion}$
is electromotive force divided by resistance. In this context, one usually
uses the concept of ``conductance'', which is the reciprocal of resistance.
The electromotive force depends on the type of charged ion (e.g., Na, K, Ca,
Cl) and the voltage of the cell; it has the form $[v-V_{ion}]$, where $%
V_{ion}$ is the so-called ``reversal potential'' of that ion.

The dynamics of the conductances are what make the equations so
mathematically interesting and rich. In a simple description of the
conductance, each ionic current has up to two ``gates'', which open or close
at rates that are dependent on the voltage of the cell. For each such gate,
there is then a first order equation of the form
\[
x^{\prime }=[x_{\infty }(v)-x]/\tau_{x}(v)
\]
where x denotes the fraction of channels of that type that are open at any
given time, $x_{\infty }(v)$ is the steady state value of $x$ for a fixed
voltage, and $\tau_{x} (v)$ is the rate constant of that gating variable. For
example, the standard Na current has the form $I_{Na}=\overline{g}%
m^{3}h\,[v-V_{Na}]$, where $\overline{g}$ is the maximal conductance, and $m$
and $h$ are gating variables satisfying the above equation for $x$. The
dynamics for $m$ and $h$ differ because $m_{\infty }(v)$ is an increasing
function of $v$, while $h_{\infty }(v)$ is a decreasing function; also $\tau
_{m}(v)$ is much smaller than $\tau _{h}(v)$. The (chemical) synaptic
currents have the same form as the intrinsic ones, with the difference that
the dependence of the driving force on voltage uses that of the
post-synaptic cell, while the conductance depends on the pre-synaptic
voltage. That is, $I_{syn}$ has the form $\overline{g}\widehat{x}%
\,[v-V_{syn}],$ where $\widehat{x}\,$satisfies the equation for x above,
with $v$ replaced by $\widehat{v}$, the voltage of the cell sending the
signal, and $V_{syn}$ is the reversal potential of the synapse. The coupling
is said to be excitatory if the current is inward (increases voltage toward
the threshold for firing an action potential) or inhibitory if the current
is outward (moves voltage away from threshold for firing.)

For a simple version of the HH equations, there are three ionic currents;
one of these (Na) creates an inward current leading to an action potential,
one (K) an outward current helping to end the action potential, and a leak
current (mainly Cl) with no gating variable. The HH equations are not one
single set of equations, but a general (and generalizable) form for a family
of equations, corresponding to different sets of intrinsic currents (which
can depend on position on the neuron), different neuron geometries, and
different networks created by interactions of neurons, which may themselves
be highly inhomogeneous. Numerical computation has become highly important
for observing the behavior of these equations, but does not suffice to
understand the behavior, especially to get insight into what the specific
ionic currents contribute; this is where the analysis, including
simplification, comes in. For an introduction to HH equations, some analysis
and some of its uses in models, see [1].

\subsection{Some mathematical issues associated with rhythms}\setzero \vskip-5mm \hspace{5mm }

It is not possible to analyze the full class of equations in all generality.
Our strategy is to look for mathematical structures underlying some classes
of behavior observed experimentally; the emphasis is on the role of
dynamical systems, as opposed to statistics, though probabilistic ideas
enter the analysis.

Our central scientific question here is how rhythms emerge from the
``wetware'', as modeled by the HH equations. As we will see, different
rhythms can be based on different sets of intrinsic currents, different
classes of neurons, and different ways of hooking up those cells. There are
some behaviors we can see by looking at small networks, and others that do
not appear until the networks are large and somewhat heterogeneous. Even in
the small networks, there are a multiplicity of different building blocks
for the rhythms, with excitation and inhibition playing different roles.
Noise appears, and plays different roles from heterogeneity.

Investigators often use simplifications of the HH equations. For example,
this talk deals only with ``space clamped'' cells in which the spatial
distribution of each cell is ignored, and the equations become ODEs. (There
are circumstances under which this can be a bad approximation, as in [2]).
Under some circumstances, the 4-dimensional simplest space-clamped HH
equations (one current equation, three gating variables) can be reduced to a
one-dimensional equation; thus, networks of neurons can be described by a
fraction of the equations that one needs for the full HH network equations.
Another kind of reduction replaces the full HH ODEs by maps that follow the
times of the spikes. In both cases, there are at least heuristic
explanations for why these reductions are often very successful, and hints
about how and why the simplifications can be expected to break down.

\section{Mathematics and small networks of neurons}\setzero \vskip-5mm \hspace{5mm}

\subsection{Centrality of inhibition in rhythms}

\vskip-5mm \hspace{5mm}

Some kinds of cells coupled by inhibition like to form rhythms and synchronize [3-5]. This is unintuitive, because
inhibition to cells can temporarily keep the latter from firing (see below for important exceptions), but mutual
inhibition can encourage cells to fire simultaneously.

There are various ways to see this, with methods that are valid in different
contexts. For weak coupling, it can be shown rigorously that the full
equations reduce to interactions between phases of the oscillators [6]; the
particular coupling associated with inhibition can then be shown to be
synchronizing (though over many cycles) [7]. If the equations can be reduced
to one-dimensional ``integrate and fire'' models, one can use
``spike-response methods'' to see the synchronizing effect of inhibitory
synapses on timing of spikes . Both of these are described in [6] along with
more references.

Another method, which I believe is most intuitive, looks at the ongoing
effect of forced inhibition on the voltage of the cells, and how some of the
processes are ``slaved'' to others. This is seen most clearly in the context
of another one-dimensional reduction that has become known as the ``theta''
model, because of the symbols used for the phase of the oscillations [8].
The reduced equations have been shown to be a canonical reduction of
equations that are near a saddle-node bifurcation on an invariant circle
(limit cycle). Many versions of HH-like models (and some kinds of real
neurons) have this property for parameter values near onset of periodic
spiking, and they are known as ``Type 1'' neurons.

The ``theta model'' has the form

\[
\theta ^{\prime }=(1-\cos \theta )+I(1+\cos \theta ).
\]
Here the equation for the phase $\theta $ has periodic solutions if the
parameter $I$ is positive, and two fixed points (stable, saddle) if $I$ is
negative. To understand the effects of forced inhibition, we replace $I$ by
a time dependent inhibition given by $I-gs(t)$, where $s(t)=\exp (-t/\tau )$
for $t>0,$ and zero otherwise. With the change of variables $J(t)=1-g\exp
(-t/\tau )$, this is a 2-D autonomous system. Figures [9] and analysis show
that the system has two special orbits, known in the non-standard analysis
literature as ``rivers'' [10], and that almost all of the trajectories feed
quickly into one of these, and are repelled from the other. The essential
effect is that initial conditions become irrelevant to the outcome of the
trajectories. A similar effect works for mutually coupled systems of
inhibitory neurons.

The rhythm formed in this way is highly dependent on the time scale of decay
of the inhibition for the frequency of the network [11, 12]. These models,
and the ``fast-firing'' inhibitory cells that they represent, can display a
large range of frequencies depending on the bias ($I_{ext}$ in HH, the
parameter $I$ in the theta model); however, in the presence of a small
amount of heterogeneity in parameters, the rhythm falls apart unless the
frequency is in the gamma range (30-80 Hz) [4, 13]. This can be understood
from spike response methods or in terms of rivers.

The above rhythm is known as ING or interneuron gamma [14, 15]. A variation
on this uses networks with fast-firing inhibitory cells (interneurons or I
cells) and excitatory cells (pyramidal cells or E-cells). This is called
PING (pyramidal interneuron gamma) [14, 15]. Heuristically, it is easy to
understand the rhythm: the inhibitory cells are set so they do not fire
without input from the E-cells. When the E-cells fire, they cause the
I-cells to cross firing threshold and inhibit the E-cells, which fire again
when the inhibition wears off. This simple mechanism becomes much more
subtle when there is heterogeneity and noise in large networks, which will
be discussed later.

\subsection{Excitation and timing maps}

\vskip-5mm \hspace{5mm}

The fast-firing cells described above are modeled using only the ionic currents needed to create a spike. Most
other neurons have channels to express many other ionic currents as well, with channel kinetics that range over a
large span of time constants. These different currents change the dynamical behavior of the cells, and allow such
cells to be ``Type II'', which means that the onset of rhythmic spiking as bias is changed is accompanied by a
Hopf bifurcation instead of a saddle node. The type of onset has important consequences for the ability of pair of
such cells to synchronize. E.g., models of the fast-firing neuron, if connected by excitatory synapses, do not
synchronize, as can be shown from weak coupling or other methods described above (e.g., [7]). However, if the
cells are Type II, they do synchronize stably with excitation (and not with inhibition). This was shown by Gutkin
and Ermentrout using weak coupling methods [16]. A more specific case study was done by Acker et al. [17],
motivated by neurons in the part of the cortex that constitutes the input-output pathways to the hippocampus, a
structure of the brain important to learning and recall. These cells are excitatory and of Type II (J.White, in
prep.); models of these cells, based on knowledge of the currents that they express, do synchronize with
excitatory synapses, and do not with inhibitory synapses.

The synchronization properties of the such cells can be understood from
spike-timing functions and maps [17]. Given the HH equations for the cell,
one can introduce at any time in the cycle excitation or inhibition whose
time course is similar to what the synapse would provide. From this, one can
compute how much the next spike is advanced or delayed by this synapse. From
such a graph, one can compute a spike-time map which takes the difference in
spike times in a single cycle to the difference in the next cycle.

The analysis of such a map is easy, but the process raises deeper
mathematical issues. One set of issues concerns what is happening at the
biophysical level that gives rise to the Type II bifurcation, which is
associated with a particular shape of the spike advance function [18].
Analysis shows that the Type II is associated with slow outward currents or certain slow inward currents that (paradoxically) turn on when the cell is inhibited [16, 17]; this shows how biophysical structure is connected with mathematical structure.

A second set of questions concerns why the high-dimensional coupled HH
equations can be well approximated by a 1-D map. (In some parameter ranges,
but not all, this is any excellent approximation). The mathematical issues
here concern how large subsets of high-dimensional phase space collapse onto
what is essentially a one-dimensional space. Ideas similar to those in
Section 3.1 are relevant, but with different biophysics creating the
collapse of the trajectories. In this case (and others) there are many
different ionic currents, with many different time scales, so that a given
current can be dominant in some portion of the trajectory and then decrease
to zero while others take over; this leads to structure that is more complex
than that of the traditional ``fast-slow'' equations, and which is not
nearly as understood. Such reductions to 1-D maps have been used in other
investigations of synchrony [19-21] involving multiple cells and multiple
kinds of currents.

\subsection{More complex building blocks: Fancier inhibitory cells}

\vskip-5mm \hspace{5mm}

So far, I've talked about networks containing fast-firing neurons (inhibitory) or excitatory cells. But there are
many different kinds of cells in the nervous system, with intrinsic and synaptic currents that make them
dynamically very different from one another. Once there are more currents with more time scales, it is easier to
create more rhythms with different frequency. That is, the differences in frequencies often (but not always) come
from some time scales in the interacting currents, and cannot be scaled away.

The stellate cell of Section 3.2 is an excellent example of currents
creating frequencies; in a wide range of parameters, these cells, even
without coupling, form a theta rhythm. Indeed, they are believed to be one
of the primary sources of that rhythm in the hippocampus, which is thought
by many to use these rhythms in tasks involving learning and recall. As
described above, these cells are excitatory, and synchronize when coupled by
excitation.

More puzzling are inhibitory cells in the hippocampus that are capable of
forming theta rhythms as isolated cells with ionic currents similar to those
in the stellate cells. The puzzle is that these cells do not cohere (in
models) using inhibitory coupling. (The decay time of inhibition caused by
these cells is roughly four times longer than the inhibition caused by the
fast-firing cells, but neither fast nor this slower decaying inhibition
creates synchrony in models.) So what is providing the coherence seen in the
theta rhythm? (The rhythm can be seen in small slices that do not have
inputs from other parts of the brain producing theta, so in such a paradigm,
the rhythm must be produced locally.)

One suggestion (Rotstein, Kopell, Whittington, in preparation) is that the
inhibitory rhythms seen in slice preparations with excitation blocked
pharmacologically depend on both kinds of inhibitory cells discussed, the
special ones (called O-LM cells [22]) and the others. Simulations show that
networks of these cells can have the O-LM cells synchronize and I-cells
synchronize at a different phase, to create an inhibitory network with
considerably more complexity than interacting fast-firing cells involved in
ING. Again, this can be reduced to a low-dimensional map for a minimal
network (two O-LM cells, one fast firing I-cell). However, the reduction now
requires properties of the currents involved in the O-LM model, including
the kinetics of the gating variables.

\subsection{Interaction of rhythms}

\vskip-5mm \hspace{5mm}

Another set of mathematical issues is associated with transitions among
rhythms. In general, rhythms slower than gamma (e.g., beta, theta and alpha)
make use of ionic currents that are active between spikes. These currents
are voltage-dependent, so that changes in voltage, in the sub- and
super-threshold regimes, can turn on or off these currents. Thus,
neuromodulators that change the voltage range of a neuron (e.g., by changing
a leak current) can change which other currents are actively expressed. In
that way, they can cause a switch from one rhythm to another. For example,
models of the alpha rhythm [20] suggest that this rhythm makes use the
inhibition-activated ``h-current''; this current is effectively off line if
the voltage is increased (even below threshold level). Thus, a switch from
alpha to a faster rhythm (gamma or beta) can be effected by simply making
the E-cells operate in a moderately higher voltage regime.

These switches can be seen in simulations (Pinto, Jones, Kaper, Kopell, in
prep.), but are still understood only heuristically. The mathematical issues
are associated with reduction of dimension methods. In the regime in which
the network is displaying alpha, there are many more variables that are
actively changing, notably the gating variables of each of the currents that
is important in this rhythm. When there is a switch to gamma and those
currents go off line, the phase space becomes effectively smaller. The
mathematics here involves understanding how that phase compression takes
place.

A related set of mathematical questions concerns rhythms that are
``nested'', one within another. For example, the theta rhythm often presents
as the envelope of a series of faster gamma cycles, and the beta rhythm, at
least in some manifestations, occur with the I-cells firing at a gamma
rhythm and the E-cells firing at the slower beta rhythm, missing some cycles
of the inhibitory rhythm. The gamma/beta switch has been understood from a
physiological point of view (see [19] and its references) and has been
simulated. The gamma/theta nesting is less understood, though new data and
simulations are providing the physiological and heuristic basis for this
[22; Rotstein, Kopell, Whittington, in prep.].
\newline

\section{Large networks}\setzero \vskip-5mm \hspace{5mm }

Though there are many more examples of other building blocks, I'm turning to
issues that do not appear in small network analysis. I'm going to go back to
a very simple building block, but now put many such together. The simple
building block is one E cell, one I-cell, which together can create a gamma
rhythm.

\subsection{Sparse coupling}

\vskip-5mm \hspace{5mm}

We now consider a network with N E-cells and M I-cells, with random coupling
from the E-cells to the I-cells and vica versa. Suppose, for example, there
is a fixed probability of connection in each direction between any pair of E
and I cells. Then the number of inputs to any cell is distributed across the
population, leading to heterogeneity of excitation and inhibition. Is it
still possible to get coherent gamma rhythms? This can be answered with
mathematical analysis using the ``theta neuron'' model described above [9].
To understand synchrony in E/I networks, it is helpful to understand what
each pulse of inhibition does to the population of excitatory cells and vica
versa. The part in which both probability and dynamical systems play a large
role is the effect of a pulse of inhibition on a population. The ``rivers''
referred to above in Section 3.1 create synchronization if the inputs to
cells have no variance, but with variation in the size of the inputs, there
is a spread in the times of the outputs. This can be accurately computed
using features of the dynamics and probability theory. Similarly, but with
less accuracy, one can compute the the effect of variation of inputs on the
spike times of the receiving population due to a pulse of excitation. The
results lead to unintutive conclusions, e.g., that increasing the strength
of the inhibition (which strengthens the synchronizing effect of the rivers)
does not reduce the desynchronizing effects of random connectivity.
Furthermore, tight synchrony can be obtained even with extremely sparse
coupling provided that variance in the size of the inputs is small.

\subsection{Loss of coherence}

\vskip-5mm \hspace{5mm}

The above analyses can be put together to understand synchrony of ``PING''.
However, they leave only partially answered many questions about larger
networks. One such question, which is central to understanding how the
assemblies of neurons are created and destroyed, is the circumstances under
which the synchrony falls apart, i.e., what modulations of cells and/or
synapses will lead to loss of coherence of the gamma rhythm. The above
analysis shows that too large a variation in size of inputs to different
cells of the same population can be fatal. Similar phenomena occur with too
much variation in drive or intrinsic currents. There are less obvious
constraints that are understood from working with smaller networks described
above. From those, it is possible to see that ING and PING operate in
different parameter regimes: the firing times of the population in ING are
governed by the bias of the I-cells (as well as the decay time of the
inhibition); in PING, the inhibitory cells are more passive until driven by
the E-cells, and the timing comes from bias of the E-cells (as well as decay
of inhibition). This means that the mechanism of coherence can switch
between ING and PING by changing relative excitability of the two
populations. Changing the strengths of the I-E and E-I synapses can also get
the population (large or small) out of the regime in which the E-cells
synchronize the I-cells, and vica versa.

A more mysterious issue that cannot be addressed within minimal networks is
how the size of the sub-populations responding on a given cycle affects the
coherence on the next cycle and the numbers of cells participating,
especially when there is some heterogeneity in the network. E.g., as the
number of inhibitory neurons firing in a cycle changes, it changes the total
inhibition to the E-cells, which changes the number of E-cells that are
ready to fire when inhibition wears off, and before the next bout of
inhibition. If the amount of inhibition gets too small, or inhibition gets
too dispersed, the coherence can rapidly die. Without taking into account
the trajectories of each of the large number of cells, it is likely that
some possibly probabilistic account of the numbers of cells spiking per cycle can
give some insight into the dynamical mechanisms surrounding failure of
coherence.

Such a reduction has been successfully used in a different setting,
involving the long-distance coherence of two populations of heterogeneous
cells. In this case, if the populations are each minimal (one E/I pair) for
each site, there is 1-D map that describes the synchronization, with the
variable the timing between the E and I sites [19] . For large and
heterogeneous networks, the synchronization (within some parameter regimes)
can be described by a 3-dimensional map, in which the first variable is the
time between the first spikes of a cycle in the two E-cell populations, and
the others are the fraction of I-cells firing on that cycle in each of the
two I-cell populations (McMillen and Kopell, in prep.).

Related work has been done from a different perspective, starting with
asynchronous networks and asking how the asynchrony can lose stability
[23-25]. Work using multiple time scales to address the formation of
``clusters'' when synchrony fails is in [26].

\subsection{Noise, PING, and frequency control}

\vskip-5mm \hspace{5mm}

One of the main differences between ING and PING is the difference in robustness. Small amounts of heterogeneity
of any kind make ING coherence fall apart dramatically [4,13]. By contrast, PING is tolerant to large ranges of
heterogeneity. The ``ping-pong'' mechanism of PING is also able to produce frequencies that cover a much wider
range than the ING mechanism, which is constrained by loss of coherence to lie in the gamma range of approximately
30-80 Hz [4,13]. Since many versions of gamma seen in experiments are of the PING variety, this raises the
question of what constrains the PING rhythms to stay in the gamma frequency range.

A possible answer to this comes from simulations. C. Borgers, D. McMillen
and I found that heterogeneity, unless extreme, would not disrupt the PING
coherence. However, a very small amount of noise (with fixed amplitude and
poisson-distributed times) could entirely destroy coherence of the PING,
provided the latter had a frequency below approximately 30 Hz; if the same
noise is introduced when the network is in the gamma range, the behavior is
only slightly perturbed. Furthermore, the ability to withstand the noise is
related to adding some I-I connections, as in ING. A heuristic explanation
is that, at low frequencies, the inhibition to the I-cells (which has a time
constant around 10ms) wears off before the excitation from the E-cells
causes these cells to spike. Thus, those cells hang around the threshold for
significant amounts of time, and are therefore vulnerable to being pushed
over threshold by noise. The mathematics has yet to be understood rigorously.

\section{Mathematics and clues to function}\setzero \vskip-5mm \hspace{5mm }

The mathematical questions are themselves interesting, but the full richness
of the scientific endeavor comes from the potential for understanding how
the rhythms generated by the brain might be used in sensory processing,
motor coordination and cognition. We are still at the outer edges of such an
investigation, but there are many clues from animal behavior, physiology and
mathematics. Work done with EEGs (see, e.g., reviews [27,28]) has shown that
many cognitive and motor tasks are associated with specific rhythms
appearing in different parts of the tasks. Gamma is often associated with
attention, awareness and perception, beta with preparation for motor
activies and high-order cognitive tasks, theta with learning and recall and
alpha with quiet awareness (there are several different versions of alpha in
different parts of the brain and found in different circumstances). Work
done in whole animals and in slice preparations are giving clues to the
underlying physiology of the rhythms, and how various neuromodulators change
the rhythms, e.g., [14]. Much of the math done so far has concerned how the
networks produce their rhythms from their ionic currents and connectivity,
and has not directly addressed function. However, the issues of function are
starting to be addressed in terms of how the dynamics of networks affects
the computational properties of the latter.

One of the potential functions for these rhythms is the creation of ``cell
assemblies'', temporary sets of neurons that fire synchronously. These
assemblies are believed to be important in distributed processing; they
enhance the effect of the synchronized pulses downstream, and provide a
substrate for changes in synapses that help to encode experience. (``Cells
that fire together wire together.'') Simulations show, and help to explain,
why gamma rhythms have especially good properties for creating cells
assemblies, and repressing cells with lower excitatbility or input [29].
Furthermore, the changes in synapses known to occur during gamma can
facilitate the creation of the beta rhythm (see [19] for references), which
appears in higher-order processing. Mathematical analysis shows that the
beta rhythm is more effective for creating synchrony over distances where
the conduction time is longer. Thus, we can understand the spontaneous
gamma-beta switch seen in various circumstances (see [19] and [29]) as
creating cell assemblies (during the gamma portion), using the synaptic
changes to get cell assemblies encoded in the beta rhythm, and then using
the beta rhythm to form highly distributed cell assemblies.

The new flood of data, plus the new insights from the mathematics, are
opening up many avenues for mathematical research related to rhythms and
function. A\ large class of such questions concerns how networks that are
displaying given rhythms filter inputs with spatio-temporal structure, and
how this affects the changing cell assemblies. This question is closely
related to the central and controversial questions of what is the neural
code and how does it operate. These questions will likely require new
techniques to combine dynamical systems and probability, new ways to reduce
huge networks to ones amenable to analysis, and new ideas within dynamical
systems itself, e.g., to understand switches as global bifurcations; these
are large and exciting challenges to the mathematical community.

\label{lastpage}

\end{document}